\documentclass[11pt,a4paper]{article}
\usepackage{authblk}
\usepackage[T1]{fontenc}
\usepackage[utf8]{inputenc}
\usepackage[british]{babel}
\usepackage{lmodern}
\usepackage[a4paper,margin=29mm]{geometry}
\usepackage{amsmath,amssymb,amsthm,mathtools}
\usepackage{booktabs,tabularx,array}
\usepackage{enumitem}
\usepackage{graphicx,float,caption}
\usepackage{microtype}
\usepackage{xcolor}
\usepackage[colorlinks=true,linkcolor=blue,citecolor=blue,urlcolor=blue]{hyperref}
\hypersetup{pdftitle={Transfer of similarity to Hilbert space contractions through inner functions: a proof of Douglas's conjecture},pdfauthor={Catalin Badea}}

\newtheorem{theorem}{Theorem}[section]
\newtheorem{proposition}[theorem]{Proposition}
\newtheorem{lemma}[theorem]{Lemma}
\newtheorem{corollary}[theorem]{Corollary}
\theoremstyle{definition}

\newtheorem{example}[theorem]{Example}

\theoremstyle{remark}
\newtheorem{remark}[theorem]{Remark}
\numberwithin{equation}{section}

\newcommand{\D}{\mathbb{D}}
\newcommand{\T}{\mathbb{T}}
\newcommand{\C}{\mathbb{C}}
\newcommand{\N}{\mathbb{N}}
\newcommand{\Hh}{\mathcal H}
\newcommand{\Bcal}{\mathcal B}
\newcommand{\simc}{\operatorname{sim}}
\newcommand{\dist}{\operatorname{dist}}
\newcommand{\supp}{\operatorname{supp}}

\newcommand{\ip}[2]{\langle #1,#2\rangle}
\newcommand{\eps}{\varepsilon}

\title{Transfer of similarity to Hilbert space contractions through inner functions: a proof of Douglas's conjecture}
\author{Catalin Badea}
\affil{%
Department of Mathematics and Statistics, University of Reading, United Kingdom,\\
and UMR 8524 -- Laboratoire Paul Painlev\'e, Universit\'e de Lille, France\\
\texttt{c.badea@reading.ac.uk} \quad  \texttt{cbadea@univ-lille.fr}\\
\url{https://www.reading.ac.uk/maths-and-stats/staff/catalin-badea}
}

\date{}

\begin{document}
\maketitle

\begin{abstract}
Let $A$ be a bounded operator on a Hilbert space and let $u$ be a
non-constant Blaschke product, finite or infinite.  We prove that if the
spectrum of $A$ does not meet the closure of the set of poles of $u$,
then $u(A)$ is similar to a contraction if and only if $A$ is.  This
settles a conjecture of R.G. Douglas. For infinite products $u$, no additional condition on the
zeros beyond the Blaschke condition is required.  We also show that $u(A)$ is similar to a contraction if and only if $A$ is
for all non-constant inner functions whose reflected singularity
set, the closure of the set of poles together with the singular
support, avoids the spectrum of $A$.
We also prove a weighted resolvent version of the theorem, which allows the poles to
accumulate on the spectrum to a limited extent; a column criterion for
partial products, which is shown to force their norm convergence; and a
simultaneous version for commuting tuples of operators tested by
separated inner functions sharing one Lyapunov operator. 
For finite Blaschke products we obtain a polynomial form of the Lyapunov operator, which gives an estimate for the similarity constants.
We also characterise the tuples of finite Blaschke products that form a uniform strong test collection: they are exactly
the tuples that are injective, with non-vanishing joint derivative, on
the closed disc.  Finally, we show by examples that the separation
hypothesis cannot be replaced by the existence of a bounded $H^\infty$
calculus, by the existence of a norm radial limit, or by the strong
convergence of the partial products, even for Ritt operators.
\end{abstract}



%%=====================================================================
\section{Introduction and preliminaries}
%%=====================================================================

\subsection{Operators similar to contractions}
Let $A\in\Bcal(\Hh)$ be a bounded linear operator on a complex Hilbert
space $\Hh$.  We say that $A$ is a \emph{contraction} if $\|A\|\le1$,
and that $A$ is \emph{similar to a contraction} if there is an
invertible operator $L\in\Bcal(\Hh)$ such that $\|L^{-1}AL\| \le 1$. 
Equivalently, $A$ is similar to a contraction if and only
if there is an equivalent Hilbertian norm on $\Hh$, that is, a norm
$\|x\|_H=\langle Hx,x\rangle^{1/2}$ with $H$ a coercive positive
operator, with respect to which $A$ is a contraction. We then call $A$
an \emph{$H$-contraction} and $H$ a \emph{Lyapunov operator} for $A$
(see Sections~\ref{sec:notation} and \ref{sec:prelim}).  By von
Neumann's
inequality, an operator similar to a contraction is \emph{polynomially
bounded}, that is,
\begin{equation}\label{eq:PB}
 \sup\{\|p(A)\| : p\in\C[z],\ \|p\|_\infty\le1\}<\infty ,
\end{equation}
where the supremum is taken over all polynomials $p$ whose supremum
norm $\|p\|_\infty:=\sup\{|p(z)|:|z|\le1\}$ is at most one.  The
estimate \eqref{eq:PB} then extends to the disc algebra.  Whether the
converse holds, that is, whether every polynomially bounded operator is
similar to a contraction, was the sixth of the ten problems posed by
Halmos \cite{HalmosTen}; Pisier \cite{Pisier1997} showed in 1997 that
the answer is negative.  On the other hand, Paulsen \cite{Pau} proved
that $A$ is similar to a contraction if and only if $A$ is
\emph{completely polynomially bounded}, which means that \eqref{eq:PB}
holds for all matrix-valued (equivalently, operator-valued) polynomials $P$,
with $\|P\|_\infty$ now denoting the supremum of the operator norms
$\|P(z)\|$ over the closed unit disc.

\subsection{Douglas's conjecture}
Let $u$ be a Blaschke product whose poles have closure disjoint from
the spectrum of $A$.  Then $u(A)$ is defined by the holomorphic
functional calculus.  If $A$ is similar to a contraction, then so is
$u(A)$: it suffices to apply von Neumann's inequality in an $H$-norm in
which $A$ is a contraction.  This paper is concerned
with the converse implication,
\begin{equation}\label{eq:Q-intro}
 u(A)\ \text{similar to a contraction}
 \quad\Longrightarrow\quad
 A\ \text{similar to a contraction},
\end{equation}
and with its analogue for more general inner functions $u$.  As
recorded in several papers, among them \cite{LanceStessin},
\cite{Stessin} and \cite{KK}, Ronald G.~Douglas conjectured in the 1990s
that the implication \eqref{eq:Q-intro} holds for every infinite
Blaschke product.  For powers, $u(z)=z^n$, the implication is
elementary and has been noted by several authors; see for instance
\cite[p.~912]{HalmosTen}.  For finite Blaschke products it was proved by
Mascioni \cite{Mascioni}, who also treated the more general setting of
Banach space operators related to contractions; see also \cite{Pis}.
Two further proofs of Mascioni's theorem in the Hilbert space case were
given by Lance and Stessin \cite{LanceStessin}, using a Wold-type
decomposition, and by Dritschel, Est\'evez and Yakubovich
\cite[Section~3.5]{DEY}.  For an infinite Blaschke product with zeros
$a_j$, Stessin \cite{Stessin} proved \eqref{eq:Q-intro} under the
additional hypothesis $\sum_j\sqrt{1-|a_j|^2}<\infty$, which is
strictly stronger than the Blaschke condition.  The hypothesis is an
artefact of the method: the Wold reconstruction series in
\cite{Stessin} is estimated by the sum of the norms of its coefficient
operators, and that sum converges only under the square-root
summability condition above.
Kazas and Kelley \cite{KK} used Clark measures to enlarge the class of
inner functions for which the transfer of similarity is known.

\subsection{Our contribution}
The aim of this paper is to prove Douglas's conjecture: the implication
\eqref{eq:Q-intro} holds for every Blaschke product $u$ and every
operator $A$ whose spectrum avoids the closure of the set of poles of
$u$, and, more generally, for every inner function $u$ whose reflected
singularity set $\operatorname{Sp}(u)$, defined in
\eqref{eq:inner-spectrum} below, avoids the spectrum of $A$.  The proof
produces, for every Lyapunov operator $X$ for $u(A)$, an explicit
coercive operator $Y$ satisfying the exact defect identity
\[
 Y-A^*YA=X-u(A)^*Xu(A),
\]
so that $Y$ is a Lyapunov operator for $A$.  Consequently, every
$H$-norm in which $u(A)$ is a contraction yields an $H$-norm in which
$A$ is one.  For Blaschke products the operator $Y$ is a
norm-convergent series of positive
operators; for the singular inner factor it is given by an integral
identity of the same kind.  We also establish three variants of the
main theorem: a weighted resolvent version, which permits the poles to
accumulate on the spectrum to a limited extent; a column criterion for
the partial products, which we show forces their norm convergence; and
a simultaneous version for commuting tuples of operators tested by
separated inner functions sharing one Lyapunov operator.  For finite
Blaschke products we derive a polynomial form of the Lyapunov operator,
leading to the estimate $\simc(A)\le C_B\simc(B(A))\max_{k<m}\|A^k\|$
for the similarity constants, and we characterise the tuples of finite
Blaschke products that form a uniform strong test collection as exactly
the tuples that are injective with non-vanishing joint derivative on
the closed disc.  To show that the separation hypothesis cannot simply
be dropped, three examples are given: a
polynomially bounded operator, not similar to a contraction, that is
annihilated by a Blaschke product; an operator that is not power
bounded although $B(rA)\to0$ in norm as $r\uparrow1$; and a Ritt
operator, not polynomially bounded, whose Blaschke partial products
converge strongly to zero.

\subsection{Organisation of the paper}
The remainder of this section fixes the notation, recalls the
holomorphic functional calculus for inner functions, and collects two
preliminary results.  Section~\ref{sec:blaschke} proves the transfer
theorem for Blaschke products, and Section~\ref{sec:inner} treats the
singular inner factor, so that every inner function whose reflected
singularity set avoids the spectrum is covered.
Section~\ref{sec:column} contains the three variants of the transfer
theorem mentioned above.  Section~\ref{sec:negative} contains the
examples showing that several natural substitutes for the separation
hypothesis are insufficient.  The last two sections concern finite
Blaschke products: Section~\ref{sec:finite} gives the polynomial form of
the Lyapunov operator, and Section~\ref{sec:uniform} the uniform strong test
theorem and its exact characterisation.  

\subsection{Notation and the functional calculus for inner functions}
\label{sec:notation}

Throughout, $\D$ is the open unit disc, $\T$ the unit circle, and
$\overline\D=\D\cup\T$.  For $a\in\D$ we put
\begin{equation}\label{eq:blaschke-factor}
 b_a(z)=\begin{cases}
 \dfrac{|a|}{a}\,\dfrac{a-z}{1-\overline az}, & a\ne0,\\[6pt]
 z, & a=0,
 \end{cases}
\end{equation}
so that $|b_a|=1$ on $\T$ and $b_a(a)=0$.  A \emph{Blaschke product} is
a product $\eta\prod_jb_{a_j}$, with $|\eta|=1$, over a finite or
infinite sequence of zeros satisfying the Blaschke condition
$\sum_j(1-|a_j|)<\infty$; it is an inner function on $\D$.  For a finite
positive singular Borel measure $\mu$ on $\T$,
\begin{equation}\label{eq:singular-factor}
 S_\mu(z)=\exp\Bigl(-\int_\T\frac{1+\overline\zeta z}{1-\overline\zeta z}
 \,d\mu(\zeta)\Bigr)
\end{equation}
is the associated singular inner function; the integral defines a
holomorphic function on $\C\setminus\supp\mu$.  Every inner function
factorises as $\theta=\eta BS_\mu$ with $|\eta|=1$, $B$ a Blaschke
product and $\mu$ singular.  Following \cite{KK}, we call
\begin{equation}\label{eq:inner-spectrum}
 \operatorname{Sp}(\theta)=\overline{\{1/\overline{a_j}:a_j\ne0\}}
 \cup\supp\mu
\end{equation}
the \emph{reflected singularity set} of $\theta$; it consists of the
poles $1/\overline{a_j}$ of $B$, their accumulation points on $\T$, and
the support of $\mu$, and $\theta$ extends holomorphically to
$\C\setminus\operatorname{Sp}(\theta)$.  The notation
$\operatorname{Sp}(\theta)$ is reserved for this set, the spectrum of an
operator $A$ being denoted by $\sigma(A)$.

We use the notation $\Bcal(\Hh)$ for the algebra of bounded linear operators on the Hilbert space
$\Hh$, and $R(\lambda,A)=(\lambda I-A)^{-1}$ for the resolvent.  The
\emph{similarity constant} of $A$ is
\begin{equation}\label{eq:sim-def}
 \simc(A)=\inf\bigl\{\|X\|\,\|X^{-1}\|:\ X\ \text{invertible},\
   \|XAX^{-1}\|\le1\bigr\},
\end{equation}
with the value $+\infty$ when $A$ is not similar to a contraction. Note
that $\sup_n\|A^n\|\le\simc(A)$.  A bounded self-adjoint operator $H$ is
\emph{coercive} if $H\ge cI$ for some $c>0$; a coercive $H$ defines the
equivalent Hilbertian norm $\|x\|_H=\langle Hx,x\rangle^{1/2}$, called
the \emph{$H$-norm}.  An operator $S$ is an \emph{$H$-contraction} if
$\|Sx\|_H\le\|x\|_H$ for all $x$, that is, if $S^*HS\le H$; in this
case we call $H$ a \emph{Lyapunov operator} for $S$. 
For a
Blaschke product $B$ with zeros $(a_j)$, listed with multiplicity, and
an operator $A$ for which the factors $I-\overline{a_j}A$ are
invertible, we write
\begin{equation}\label{eq:PE-notation}
 P_0=I,\qquad P_N=\prod_{j=1}^Nb_{a_j}(A),\qquad
 E_j=\sqrt{1-|a_j|^2}\,(I-\overline{a_j}A)^{-1}P_{j-1}
\end{equation}
for the partial products of $B$ at $A$ and for the factors of the
Lyapunov operators constructed below.

We recall how $\theta(A)$ is defined for an inner function $\theta$
whose reflected singularity set avoids $\sigma(A)$; compare
\cite[Section~5]{LanceStessin}.  Let $f$ be holomorphic on an open
neighbourhood $U$ of $\sigma(A)$, and let $V$ be an open neighbourhood
of $\sigma(A)$ with compact closure in $U$ whose boundary
$\Gamma=\partial V$ is a finite union of smooth Jordan curves.  The
Riesz--Dunford functional calculus defines
\begin{equation}\label{eq:riesz-dunford}
 f(A)=\frac1{2\pi i}\int_\Gamma f(z)\,R(z,A)\,dz .
\end{equation}
Let $\theta=\eta BS_\mu$ be an inner function and suppose that
\begin{equation}\label{eq:separation-intro}
 \sigma(A)\cap\operatorname{Sp}(\theta)=\varnothing .
\end{equation}
The Blaschke product converges locally uniformly on
$\C\setminus\operatorname{Sp}(B)$, where it is holomorphic, and
\eqref{eq:singular-factor} is holomorphic on $\C\setminus\supp\mu$, so
$\theta$ is holomorphic on the open set
$\C\setminus\operatorname{Sp}(\theta)\supset\sigma(A)$ and
\eqref{eq:riesz-dunford} defines $\theta(A)=\eta B(A)S_\mu(A)$.  It is
easy to see that the product
\begin{equation}\label{eq:B-of-A-product}
 B(A)=\eta\prod_{j}b_{a_j}(A),\qquad
 b_{a_j}(A)=\frac{|a_j|}{a_j}\,(a_j I-A)(I-\overline{a_j}A)^{-1}
 \quad(a_j\ne0),
\end{equation}
converges in operator norm; this is proved in
Lemma~\ref{lem:convergence-p2} below.  With this notation, the problem
raised by Douglas reads as follows: \emph{let $\theta$ be an inner
function satisfying \eqref{eq:separation-intro}; if $\theta(A)$ is
similar to a contraction, is $A$ similar to a contraction?}  This is
the form in which the problem is stated in
\cite[Section~5]{LanceStessin}; Theorem~\ref{thm:inner} answers it in
the affirmative.

Finally, for a finite non-empty set $E\subset\T$ we call $A$ a
$\mathrm{Ritt}_E$ operator if $\sigma(A)\subset\overline\D$ and
$\|R(\lambda,A)\|\le C/\dist(\lambda,E)$ for some $C<\infty$ and all
$|\lambda|>1$.  A \emph{Ritt operator} is an operator with
$\sigma(A)\subset\overline\D$ and $\|R(\lambda,A)\|\le C/|\lambda-1|$
for $|\lambda|>1$; equivalently \cite{NZ,Lyu}, $A$ is power bounded and
$\sup_kk\|A^{k-1}(I-A)\|<\infty$.

\subsection{Two preliminary results}\label{sec:prelim}

The first result is the standard reformulation of similarity to a
contraction in terms of Lyapunov operators; it is the form in which
similarity will be established throughout the paper.

\begin{proposition}[Lyapunov criterion]\label{prop:metric-p2}
For $S\in\Bcal(\Hh)$ the following are equivalent:
\begin{enumerate}[label=\textup{(\roman*)}]
\item $S$ is similar to a contraction;
\item $S$ is an $H$-contraction for some coercive $H$, that is, $S$
      admits a Lyapunov operator.
\end{enumerate}
If $H\ge cI$ with $c>0$ and $S^*HS\le H$, then $H^{1/2}SH^{-1/2}$ is a
contraction and $\simc(S)\le\|H\|^{1/2}c^{-1/2}$; conversely, if
$XSX^{-1}$ is a contraction, then $H=X^*X$ is a Lyapunov operator for
$S$.  Moreover, if $\simc(S)<\infty$, then the infimum in
\eqref{eq:sim-def} is attained.
\end{proposition}

\begin{proof}
The equivalence follows by conjugating the inequality $S^*HS\le H$ by
$H^{-1/2}$.  For the attainment of the infimum, let $K=\simc(S)$ and
choose Lyapunov operators $H_n$ for $S$, scaled so that
$I\le H_n\le(K+1/n)^2I$.  Any ultraweak cluster point $H$ of the
bounded sequence $(H_n)$ satisfies $I\le H\le K^2I$ and $S^*HS\le H$,
because these inequalities define ultraweakly closed sets, and the
similarity $H^{1/2}SH^{-1/2}$ has condition number at most $K$.
\end{proof}

Everything in Sections~\ref{sec:blaschke}--\ref{sec:column} rests on the
following consequence of the separation hypothesis.  Note that the lemma
does \emph{not} assume $\sigma(A)\subset\overline\D$: in the theorems
below this inclusion is a consequence of the conclusion.

\begin{lemma}[Norm convergence and uniform bounds]\label{lem:convergence-p2}
Let $B=\eta\prod_{j\in J}b_{a_j}$ be a non-constant Blaschke product
with $J=\{1,\dots,m\}$ or $J=\N$, its zeros being listed with
multiplicity and the finitely many zeros at the origin coming first, and
suppose that
\begin{equation}\label{eq:separation}
 \sigma(A)\cap\operatorname{Sp}(B)=\varnothing .
\end{equation}
Then every $I-\overline{a_j}A$ is invertible and, with $P_N$ as in
\eqref{eq:PE-notation},
\begin{equation}\label{eq:PLbeta}
 L:=\sup_j\|(I-\overline{a_j}A)^{-1}\|<\infty,\quad
 P:=\sup_{N}\|P_N\|<\infty,\quad
 \beta:=\sum_j(1-|a_j|^2)<\infty .
\end{equation}
In the infinite case the sequence $(P_N)$ converges in operator norm,
so that $B(A)=\eta\lim_NP_N$ is well defined, and this limit coincides
with the operator given by
the functional calculus \eqref{eq:riesz-dunford}.
\end{lemma}

\begin{proof}
The set $\overline{\{a_j\}}\subset\overline\D$ is compact, and for every
$a$ in it the operator $I-\overline aA$ is invertible: for $a=0$ this is
clear, and for $a\ne0$ the point $1/\overline a$ lies in
$\operatorname{Sp}(B)$, hence outside $\sigma(A)$.  Since inversion is
continuous on the group of invertible operators, $L<\infty$.

For $a=re^{i\theta}\ne0$ one has the identity
\begin{equation}\label{eq:factor-close}
 1-b_a(z)=(1-r)\,(1+e^{-i\theta}z)\,(1-\overline az)^{-1},
\end{equation}
which is obtained by writing the left-hand side over the common
denominator $1-\overline az$ and using $\overline a=re^{-i\theta}$.
Hence
\[
 \sum_{j:a_j\ne0}\|I-b_{a_j}(A)\|\le(1+\|A\|)L\sum_{j:a_j\ne0}(1-|a_j|)
 <\infty .
\]
The zeros at the origin contribute a finite product.  If $F_j$ denotes
the remaining factors, then $\sum_j\|F_j-I\|<\infty$ and
$\|\prod_{j=N+1}^MF_j-I\|\le\exp\bigl(\sum_{j=N+1}^M\|F_j-I\|\bigr)-1$,
so the partial products form a Cauchy sequence and are bounded.  Writing
$q$ for the number of zeros at the origin, we obtain
\[
 P\le\max\{1,\|A\|\}^q\exp\Bigl((1+\|A\|)L\sum_{j:a_j\ne0}(1-|a_j|)\Bigr).
\]
Finally, $1-|a_j|^2\le2(1-|a_j|)$, so $\beta<\infty$.  The same
estimates hold uniformly on compact subsets of the complement of
$\operatorname{Sp}(B)$, so the scalar product converges locally
uniformly to a holomorphic function there, and evaluating it at $A$ by
means of a fixed contour around $\sigma(A)$ shows that the norm limit of
$(P_N)$ agrees with the functional calculus.
\end{proof}

%%=====================================================================
\section{The transfer theorem for Blaschke products}
\label{sec:blaschke}
%%=====================================================================

Throughout this section $B=\eta\prod_{j\in J}b_{a_j}$ is a non-constant
Blaschke product, finite or infinite, and $A$ satisfies the separation
hypothesis \eqref{eq:separation}.  As in \eqref{eq:PE-notation}, we put
\begin{equation}\label{eq:E-def}
 E_j=\sqrt{1-|a_j|^2}\,(I-\overline{a_j}A)^{-1}P_{j-1},
 \qquad j\in J,
\end{equation}
where $P_{j-1}=\prod_{k<j}b_{a_k}(A)$ and $P_0=I$.  All the operators
$P_j$ and $E_j$ are rational functions of $A$, and therefore commute
with $A$ and with one another.

\begin{theorem}[Transfer for every Blaschke product]\label{thm:blaschke}
Let $B$ and $A$ be as above, so that \eqref{eq:separation} holds.  Then
\begin{equation}\label{eq:blaschke-iff}
 B(A)\ \text{is similar to a contraction}
 \iff A\ \text{is similar to a contraction}.
\end{equation}
More precisely, if $H$ is a Lyapunov operator for $B(A)$, then the
series
\begin{equation}\label{eq:Q-series}
 Q=\sum_{j\in J}E_j^*HE_j
\end{equation}
converges in operator norm, its sum is a Lyapunov operator for $A$, and
\begin{equation}\label{eq:stein-blaschke}
 Q-A^*QA=H-B(A)^*HB(A).
\end{equation}
\end{theorem}

\begin{proof}
\emph{Convergence and coercivity.}  By Lemma~\ref{lem:convergence-p2},
\begin{equation}\label{eq:E-sum}
 \sum_{j\in J}\|E_j\|^2\le P^2L^2\sum_j(1-|a_j|^2)=P^2L^2\beta<\infty ,
\end{equation}
so the series of positive operators \eqref{eq:Q-series} converges in
norm, with the tail estimate
$\|Q-Q_N\|\le\|H\|P^2L^2\sum_{j>N}(1-|a_j|^2)$ for the partial sums
$Q_N=\sum_{j\le N}E_j^*HE_j$, and
\begin{equation}\label{eq:Q-upper}
 0\le Q\le\|H\|P^2L^2\beta\,I .
\end{equation}
The first term of the series already gives coercivity, since
$E_1=\sqrt{1-|a_1|^2}\,(I-\overline{a_1}A)^{-1}$ is invertible:
\begin{equation}\label{eq:Q-lower}
 Q\ge E_1^*HE_1\ge
 \frac{1-|a_1|^2}{\|H^{-1}\|\,\|I-\overline{a_1}A\|^2}\,I>0 .
\end{equation}

\emph{The identity.}  Let $|a|<1$ and let $H$ be self-adjoint.  Put
$c=I-\overline aA$ and $b(z)=(z-a)/(1-\overline az)$, a unimodular
multiple of $b_a$.  Since $b(A)c=c\,b(A)=A-aI$, all the factors being
rational functions of $A$,
\begin{equation}\label{eq:one-factor-p2}
\begin{split}
 c^*(H-b(A)^*Hb(A))c
 &=c^*Hc-(b(A)c)^*H(b(A)c)\\
 &=\bigl(I-aA^*\bigr)H\bigl(I-\overline aA\bigr)
  -(A^*-\overline aI)H(A-aI)\\
 &=(1-|a|^2)(H-A^*HA).
\end{split}
\end{equation}
Since a unimodular factor does not affect $b_a(A)^*Hb_a(A)$, this can be
written as
\begin{equation}\label{eq:one-factor-b}
 H-b_a(A)^*Hb_a(A)=(1-|a|^2)(I-\overline aA)^{-*}
   (H-A^*HA)(I-\overline aA)^{-1}.
\end{equation}
Since $P_j=b_{a_j}(A)P_{j-1}$ and all the factors commute, telescoping
gives, for every $N$,
\[
 H-P_N^*HP_N=\sum_{j=1}^NP_{j-1}^*\bigl(H-b_{a_j}(A)^*Hb_{a_j}(A)\bigr)P_{j-1},
\]
and inserting \eqref{eq:one-factor-b} we obtain
\begin{equation}\label{eq:finite-telescoping}
 H-P_N^*HP_N=\sum_{j=1}^NE_j^*(H-A^*HA)E_j=Q_N-A^*Q_NA ,
\end{equation}
the last equality because each $E_j$ commutes with $A$.  Both sides
converge in norm as $N\to\infty$: $P_N\to\overline\eta B(A)$ by
Lemma~\ref{lem:convergence-p2}, and $Q_N\to Q$ by \eqref{eq:E-sum}.
This proves \eqref{eq:stein-blaschke}.

\emph{Conclusion.}  The right-hand side of \eqref{eq:stein-blaschke} is
non-negative by hypothesis, so $A^*QA\le Q$, and $Q^{1/2}AQ^{-1/2}$ is
a contraction by Proposition~\ref{prop:metric-p2}.

\emph{Converse.}  If $XAX^{-1}=C$ is a contraction, then
$P_N(C)=XP_N(A)X^{-1}$ is a contraction for every $N$ by von Neumann's
inequality, and $\overline\eta XB(A)X^{-1}=\lim_NP_N(C)$ in norm by
Lemma~\ref{lem:convergence-p2}; hence $XB(A)X^{-1}$ is a contraction.
\end{proof}

\begin{remark}[Why square-root summability is not needed]
Estimating the terms of \eqref{eq:E-def} individually gives
$\|E_j\|\le PL\sqrt{1-|a_j|^2}$, and a bound for $\sum_j\|E_j\|$
requires $\sum_j\sqrt{1-|a_j|^2}<\infty$.  The series
\eqref{eq:Q-series}, however, is a series of \emph{positive} operators,
whose norm is controlled by $\sum_j\|E_j\|^2$, and this costs only
$\sum_j(1-|a_j|^2)<\infty$, which is the Blaschke condition.  The
passage from this estimate to the operator identity is the finite
telescoping identity \eqref{eq:finite-telescoping}, which involves no
infinite sums at all.  This is what removes the square-root
summability hypothesis of \cite{Stessin}.
\end{remark}

\begin{corollary}[Quantitative similarity bound]\label{cor:constant-p2}
In the notation of Theorem~\ref{thm:blaschke}, the similarity
$Q^{1/2}$ satisfies
\begin{equation}\label{eq:cond-number}
 \|Q^{1/2}\|\,\|Q^{-1/2}\|\le
 \sqrt{\|H\|\,\|H^{-1}\|}\,PL\,\|I-\overline{a_1}A\|
 \sqrt{\frac{\beta}{1-|a_1|^2}} .
\end{equation}
If $a_1=0$, then $E_1=I$, so $Q\ge H$, and the bound reads
$\sqrt{\|H\|\|H^{-1}\|}\,PL\sqrt\beta$.
\end{corollary}

\begin{proof}
Combine \eqref{eq:Q-upper} with \eqref{eq:Q-lower}; when $a_1=0$ one
has $E_1=P_0=I$.
\end{proof}

\begin{remark}\label{rem:a1}
The lower bound \eqref{eq:Q-lower} depends on the first zero only,
whereas $P=\sup_N\|P_N\|$ may depend on the order in which the zeros
are listed; the exponential bound for $P$ in
Lemma~\ref{lem:convergence-p2}, however, does not.  Listing a zero at
the origin first removes the factor $\|I-\overline{a_1}A\|$ from
\eqref{eq:cond-number}.  The constants $P$ and $L$ depend on $A$.  For
finite products, a bound in terms of $\simc(B(A))$ and finitely many
powers of $A$ is given in Section~\ref{sec:finite}, and the dependence
on the powers cannot be removed (Remark~\ref{rem:exponent}).
\end{remark}

\begin{remark}[On the hypotheses]
The assumption that $B$ is non-constant cannot be dropped: a unimodular
constant $B$ gives $B(A)=\eta I$ for every $A$.  The separation
hypothesis \eqref{eq:separation} cannot be replaced by the mere
invertibility of the factors $I-\overline{a_j}A$, as
Section~\ref{sec:negative} shows; it can, however, be weakened to a
weighted resolvent condition, see Section~\ref{subsec:weighted}.
\end{remark}

%%=====================================================================
\section{Singular inner factors}
\label{sec:inner}
%%=====================================================================

Let $\mu$ be a finite positive singular Borel measure on $\T$ and let
$S_\mu$ be the singular inner function \eqref{eq:singular-factor}.  The
separation hypothesis now reads
\begin{equation}\label{eq:singular-separation}
 \sigma(A)\cap\supp\mu=\varnothing .
\end{equation}
For $\zeta\in\supp\mu$ put $D_\zeta=(I-\overline\zeta A)^{-1}$, and
define
\begin{equation}\label{eq:Fmu}
 F=F_\mu(A)=\int_\T(I+\overline\zeta A)D_\zeta\,d\mu(\zeta),
 \qquad E_s=e^{-sF}\quad(0\le s\le1),
\end{equation}
all integrals being Bochner integrals in the operator norm.  By
\eqref{eq:singular-separation} and the compactness of $\supp\mu$,
\begin{equation}\label{eq:Lmu}
 L_\mu:=\sup_{\zeta\in\supp\mu}\|D_\zeta\|<\infty
\end{equation}
when $\mu\ne0$, so $F$ and $E_1=e^{-F}$ are well defined.  To identify
$E_1$ with the operator given by the functional calculus, put
$\Lambda(z)=-\int_\T\frac{1+\overline\zeta z}{1-\overline\zeta z}\,d\mu(\zeta)$,
which is holomorphic on $\C\setminus\supp\mu$ and satisfies
$S_\mu=e^{\Lambda}$ there.  Since $\sigma(A)$ is a compact subset of
$\C\setminus\supp\mu$ and the functional calculus commutes with
norm-convergent Bochner integrals, $\Lambda(A)=-F$, and hence
$S_\mu(A)=e^{-F}=E_1$.  When $\mu=0$ we set $F=0$ and $S_\mu(A)=I$.

\begin{lemma}[A Lyapunov operator for the singular factor]\label{lem:singular}
Assume \eqref{eq:singular-separation} and let $H=H^*$ be bounded.  Put
\begin{equation}\label{eq:Qmu}
 Q_\mu(H)=2\int_0^1\int_\T D_\zeta^*E_s^*HE_sD_\zeta
 \,d\mu(\zeta)\,ds .
\end{equation}
Then
\begin{equation}\label{eq:stein-singular}
 Q_\mu(H)-A^*Q_\mu(H)A=H-S_\mu(A)^*HS_\mu(A).
\end{equation}
If $\mu\ne0$ and $H\ge hI$ with $h>0$, then $Q_\mu(H)$ is coercive;
more precisely, with $t=\mu(\T)$ and $f=\|F\|$,
\begin{equation}\label{eq:Qmu-bounds}
 \frac{2the^{-2f}}{(1+\|A\|)^2}\,I\le Q_\mu(H)
 \le2tL_\mu^2e^{2f}\|H\|\,I .
\end{equation}
\end{lemma}

\begin{proof}
For $\zeta\in\supp\mu$ put $F_\zeta=(I+\overline\zeta A)D_\zeta$.
Multiplying $F_\zeta^*H+HF_\zeta$ on the left by $I-\zeta A^*$ and on
the right by $I-\overline\zeta A$ gives, since $|\zeta|=1$,
\[
 (I+\zeta A^*)H(I-\overline\zeta A)+(I-\zeta A^*)H(I+\overline\zeta A)
 =2(H-A^*HA),
\]
that is,
\begin{equation}\label{eq:herglotz}
 F_\zeta^*H+HF_\zeta=2D_\zeta^*(H-A^*HA)D_\zeta .
\end{equation}
Integrating \eqref{eq:herglotz} against $\mu$ and using the identity
$\frac{d}{ds}(E_s^*HE_s)=-E_s^*(F^*H+HF)E_s$, which follows from
$E_s'=-FE_s$ and the fact that $F$ commutes with $E_s$, we obtain
\[
 H-E_1^*HE_1=2\int_0^1\int_\T(E_sD_\zeta)^*(H-A^*HA)(E_sD_\zeta)
 \,d\mu(\zeta)\,ds .
\]
All the factors are functions of $A$ and hence commute with $A$, so the
right-hand side equals $Q_\mu(H)-A^*Q_\mu(H)A$.

As for the bounds, the integrand in \eqref{eq:Qmu} is jointly norm
continuous on $[0,1]\times\supp\mu$ and bounded by
$L_\mu^2e^{2f}\|H\|$, which gives the existence of the integral and the
upper estimate; positivity is clear from the definition.  Each
$E_sD_\zeta$ is invertible, with
$\|(E_sD_\zeta)^{-1}\|=\|(I-\overline\zeta A)e^{sF}\|\le(1+\|A\|)e^f$,
so that
$\ip{D_\zeta^*E_s^*HE_sD_\zeta x}{x}\ge
h\,e^{-2f}(1+\|A\|)^{-2}\|x\|^2$, and integrating gives the lower
bound.
\end{proof}

\begin{theorem}[Transfer for separated inner functions]\label{thm:inner}
Let $\theta=\eta BS_\mu$ be a non-constant inner function, where $B$ is
a finite or infinite Blaschke product or the empty product and $\mu$ is
a finite positive singular measure, and suppose that
\begin{equation}\label{eq:inner-separation}
 \sigma(A)\cap\operatorname{Sp}(\theta)=\varnothing .
\end{equation}
Then $\theta(A)$ is well defined, and
\begin{equation}\label{eq:inner-iff}
 \theta(A)\ \text{is similar to a contraction}
 \iff A\ \text{is similar to a contraction}.
\end{equation}
If $H$ is a Lyapunov operator for $\theta(A)$, then
\begin{equation}\label{eq:combined-Q}
 Q=Q_B(H)+B(A)^*Q_\mu(H)B(A),
\end{equation}
where $Q_B(H)$ is the series \eqref{eq:Q-series} and $Q_B(H)=0$ if $B$
is the empty product, is a Lyapunov operator for $A$, and
\begin{equation}\label{eq:stein-inner}
 Q-A^*QA=H-\theta(A)^*H\theta(A).
\end{equation}
\end{theorem}

\begin{proof}
Both constructions apply under \eqref{eq:inner-separation}, since
$\operatorname{Sp}(\theta)$ contains the closure of the set of poles of
$B$ as well as $\supp\mu$; the series and the integral converge, their
terms are positive, and $\theta(A)=\eta B(A)S_\mu(A)$.  If $B$ is
non-constant, then $Q_B(H)$ is coercive by \eqref{eq:Q-lower}; if $B$
is the empty product, then $\mu\ne0$ because $\theta$ is non-constant,
and \eqref{eq:Qmu-bounds} gives the coercivity of $Q_\mu(H)$.  In either
case $Q$ is coercive.

Since $A$ commutes with $B(A)$,
\begin{align*}
 Q-A^*QA&=\bigl(H-B(A)^*HB(A)\bigr)
 +B(A)^*\bigl(H-S_\mu(A)^*HS_\mu(A)\bigr)B(A)\\
 &=H-\theta(A)^*H\theta(A),
\end{align*}
the unimodular constant $\eta$ playing no role.  Hence $A^*QA\le Q$,
and $A$ is similar to a contraction by Proposition~\ref{prop:metric-p2}.

Conversely, let $RAR^{-1}=C$ be a contraction.  Then
$B(C)=RB(A)R^{-1}$ is a contraction by the converse part of
Theorem~\ref{thm:blaschke}, and \eqref{eq:herglotz} with $H=I$ gives
\[
 F_\mu(C)^*+F_\mu(C)=2\int_\T(I-\overline\zeta C)^{-*}
 (I-C^*C)(I-\overline\zeta C)^{-1}\,d\mu(\zeta)\ge0,
\]
whence $\frac{d}{ds}\bigl(e^{-sF_\mu(C)}\bigr)^*e^{-sF_\mu(C)}\le0$ and
$\|e^{-F_\mu(C)}\|\le1$.  Therefore $\theta(C)=R\theta(A)R^{-1}$ is a
contraction.
\end{proof}

\begin{example}[An atomic singular factor]\label{ex:atomic}
For $t>0$ and $\zeta\in\T$, the inner function
$\theta(z)=\exp\bigl(-t(1+\overline\zeta z)/(1-\overline\zeta z)\bigr)$
has $\operatorname{Sp}(\theta)=\{\zeta\}$, so the only hypothesis is
$\zeta\notin\sigma(A)$.  Thus
\[
 \exp\bigl(-t(I+\overline\zeta A)(I-\overline\zeta A)^{-1}\bigr)
\]
is similar to a contraction if and only if $A$ is, and an explicit
Lyapunov operator for $A$ is
\[
 Q=2t\int_0^1D_\zeta^*e^{-sF^*}He^{-sF}D_\zeta\,ds,\qquad
 F=t(I+\overline\zeta A)D_\zeta .
\]
\end{example}

\begin{remark}
Theorem~\ref{thm:inner} covers arbitrary Blaschke sequences and
arbitrary singular measures, but only when the spectrum avoids the
reflected singularity set \eqref{eq:inner-separation}.  It says nothing
about an $H^\infty$ functional calculus at a boundary singularity of
$\theta$, and Section~\ref{sec:negative} shows that no such statement
can hold in general.  Only the positivity and finiteness of $\mu$ are
used in the proof; singularity of $\mu$ serves only to identify
$S_\mu$ as an inner function.
\end{remark}

%%=====================================================================
\section{Three variants of the transfer theorem for Blaschke products}
\label{sec:column}
%%=====================================================================

The proof of Theorem~\ref{thm:blaschke} uses less than separation of
the spectrum from the whole closed set of poles: it needs a meaning for
$B(A)$ and a summability property of the factors $E_j$.  In this
section we vary these two ingredients and obtain three variants of the
theorem: a weighted resolvent hypothesis under which the proof goes
through verbatim; a column criterion for the partial products, which
turns out to force their norm convergence; and a simultaneous version
for commuting tuples of operators and arbitrary separated inner
functions.

\subsection{A weighted resolvent hypothesis}
\label{subsec:weighted}

What the proof of Theorem~\ref{thm:blaschke} requires is the norm
convergence of the partial products and the square summability of the
factors $E_j$.  Both follow from a single weighted condition on the
resolvents.

\begin{proposition}\label{prop:weighted}
Let $B$ be a non-constant Blaschke product with zeros $(a_j)_{j\in J}$,
and suppose that every $I-\overline{a_j}A$ is invertible and that
\begin{equation}\label{eq:weighted-hyp}
 \mathcal E_A(B):=\sum_{j\in J}(1-|a_j|^2)\|(I-\overline{a_j}A)^{-1}\|^2
 <\infty .
\end{equation}
Then the partial products $P_N$ converge in operator norm, so that
$B(A)=\eta\lim_NP_N$ is defined; the operator $Q$ of
\eqref{eq:Q-series} is well defined and coercive, and
\eqref{eq:stein-blaschke} holds.  Consequently \eqref{eq:blaschke-iff}
holds, and
\begin{equation}\label{eq:weighted-bounds}
 \frac{1-|a_1|^2}{\|H^{-1}\|\,\|I-\overline{a_1}A\|^2}\,I\le Q
 \le\|H\|P^2\mathcal E_A(B)\,I ,
\end{equation}
so that $L\sqrt\beta$ may be replaced by $\sqrt{\mathcal E_A(B)}$ in
\eqref{eq:cond-number}.
\end{proposition}

\begin{proof}
Let $d=\sum_{j:a_j\ne0}(1-|a_j|)<\infty$.  By the Cauchy--Schwarz
inequality and $(1-|a_j|)/(1+|a_j|)\le1-|a_j|$,
\begin{align*}
 \sum_{j:a_j\ne0}&(1-|a_j|)\|(I-\overline{a_j}A)^{-1}\|\\
 &\le\Bigl(\sum_{j:a_j\ne0}\frac{1-|a_j|}{1+|a_j|}\Bigr)^{1/2}
    \Bigl(\sum_{j:a_j\ne0}(1-|a_j|^2)
    \|(I-\overline{a_j}A)^{-1}\|^2\Bigr)^{1/2}
 \le\sqrt{d\,\mathcal E_A(B)} .
\end{align*}
By \eqref{eq:factor-close}, the product therefore converges absolutely
in norm after its finite initial part, and
$P\le\max\{1,\|A\|\}^q\exp\bigl((1+\|A\|)\sqrt{d\mathcal E_A(B)}\bigr)
<\infty$.  Moreover $\sum_j\|E_j\|^2\le P^2\mathcal E_A(B)<\infty$, so
the remainder of the proof of Theorem~\ref{thm:blaschke} applies
without change.  In the converse direction, $P_N(C)=XP_N(A)X^{-1}$ is a
contraction for every $N$ by von Neumann's inequality and converges in
norm, so no separation hypothesis is needed there either.
\end{proof}

\begin{example}[The weighted condition is strictly weaker]\label{ex:weighted}
Let $A=I$ and $a_j=(1-2^{-3j})e^{i2^{-j}}$.  These zeros satisfy the
Blaschke condition, and the poles $1/\overline{a_j}$ accumulate at
$1\in\sigma(A)$, so \eqref{eq:separation} fails.  Writing
$r_j=1-2^{-3j}$ and $\theta_j=2^{-j}$, we have
\[
 |1-a_j|^2=(1-r_j)^2+4r_j\sin^2(\theta_j/2)\ge c\,2^{-2j},
 \qquad 1-|a_j|^2\le2\cdot2^{-3j},
\]
with an absolute constant $c>0$, so that
$\mathcal E_I(B)\le C\sum_j2^{-3j}/2^{-2j}=C\sum_j2^{-j}<\infty$.
Proposition~\ref{prop:weighted} therefore applies to a Blaschke product
whose poles accumulate on the spectrum.  We do not claim that the
weighted condition is necessary.
\end{example}

\subsection{The column criterion forces norm convergence}
\label{subsec:column}

The second variant replaces the summability of the norms $\|E_j\|^2$
by a weaker, columnwise, condition.  This condition turns out to be
considerably more useful than a mere bound for strongly convergent
products: combined with the Blaschke condition, it yields the norm
convergence of the partial products, with an explicit rate.

\begin{proposition}[Column criterion]\label{prop:column}
Let $(a_j)_{j\ge1}\subset\D$ satisfy the Blaschke condition, with the
finitely many zeros at the origin listed first, and suppose that every
$I-\overline{a_j}A$ is invertible.  If
\begin{equation}\label{eq:column-hyp}
 C:=\sup_N\Bigl\|\sum_{j=1}^NE_j^*E_j\Bigr\|<\infty,
\end{equation}
then $P_N$ converges in operator norm to an operator $S$, and for every
$N$ beyond the zeros at the origin
\begin{equation}\label{eq:column-tail}
 \|S-P_N\|\le(1+\|A\|)\sqrt{C}
 \Bigl(\sum_{j>N}\frac{1-|a_j|}{1+|a_j|}\Bigr)^{1/2}.
\end{equation}
If $S$ is similar to a contraction, then so is $A$.  More precisely,
whenever $H$ is a Lyapunov operator for $S$, the strong limit
$Q=\operatorname{s-lim}_N\sum_{j\le N}E_j^*HE_j$ exists, is a Lyapunov
operator for $A$, and satisfies
\begin{equation}\label{eq:column-stein}
 Q-A^*QA=H-S^*HS\ge0.
\end{equation}
Conversely, if $A$ is similar to a contraction, then so is $S$.
\end{proposition}

\begin{proof}
For $a_j=r_je^{i\theta_j}\ne0$, formula~\eqref{eq:factor-close} gives
\[
 P_{j-1}-P_j=\alpha_j(I+e^{-i\theta_j}A)E_j,
 \qquad \alpha_j=\sqrt{\frac{1-r_j}{1+r_j}}.
\]
Hence, for $M>N$ beyond the zeros at the origin and $x\in\Hh$, the
Cauchy--Schwarz inequality gives
\begin{align*}
 \|(P_M-P_N)x\|
 &\le\Bigl(\sum_{j=N+1}^M\alpha_j^2\Bigr)^{1/2}
       \Bigl(\sum_{j=N+1}^M
       \|(I+e^{-i\theta_j}A)E_jx\|^2\Bigr)^{1/2}\\
 &\le(1+\|A\|)\sqrt C\,\|x\|
       \Bigl(\sum_{j>N}\frac{1-r_j}{1+r_j}\Bigr)^{1/2}.
\end{align*}
The last sum tends to zero by the Blaschke condition, which proves the
norm convergence of $(P_N)$ and the estimate~\eqref{eq:column-tail}.

The positive operators $Q_N=\sum_{j\le N}E_j^*HE_j$ increase with $N$
and satisfy $0\le Q_N\le\|H\|CI$, so they converge strongly to a
bounded positive operator $Q$, and $Q\ge E_1^*HE_1$ is coercive
by~\eqref{eq:Q-lower}.  The finite identity
\eqref{eq:finite-telescoping}, $Q_N-A^*Q_NA=H-P_N^*HP_N$, holds for
every $N$ without any convergence hypothesis; passing to the limit,
strongly on the left-hand side and in norm on the right-hand side,
proves~\eqref{eq:column-stein}.  Proposition~\ref{prop:metric-p2}
completes the proof of the forward implication.  Conversely, in an
$H$-norm in which $A$ is a contraction all finite products $P_N$ are
contractions by von Neumann's inequality, and so is their norm limit.
\end{proof}

\begin{remark}\label{rem:column-mean}
Condition~\eqref{eq:column-hyp} says that $x\mapsto(E_jx)_{j\ge1}$ is a
bounded map from $\Hh$ to $\ell^2(\Hh)$.  It is weaker than
$\sum_j\|E_j\|^2<\infty$, but strictly stronger than the norm
convergence of the partial products.  To see this, take $A=I$ on a
one-dimensional space and $a_j=(1-\varepsilon_j^2)e^{i\varepsilon_j}$,
where $0<\varepsilon_j<1/2$, $\varepsilon_j\downarrow0$ and
$\sum_j\varepsilon_j<\infty$.  Then $|P_N|=1$ for every $N$, while
\eqref{eq:factor-close} gives
\[
 |1-b_{a_j}(1)|
 =\frac{2\varepsilon_j^2\cos(\varepsilon_j/2)}{|1-a_j|}
 \le C_0\varepsilon_j,
\]
because $|1-a_j|\ge c\varepsilon_j$; thus $(P_N)$ converges.  On the
other hand $|1-a_j|^2\le2\varepsilon_j^2$ and
$1-|a_j|^2\ge\varepsilon_j^2$, so
$\|E_j\|^2=(1-|a_j|^2)/|1-a_j|^2\ge1/2$, and the column series
diverges.  In this example $A$ and the limit $S$ are unimodular
scalars, so the failure of the column estimate does not by itself
obstruct the transfer of similarity.
\end{remark}

The case in which the norm limit has spectral radius less than one is
elementary and requires neither separation nor the column estimate.

\begin{proposition}[A strictly stable norm limit]\label{prop:stable-limit}
Suppose that every $I-\overline{a_j}A$ is invertible and that
$P_N\to S$ in operator norm.  If the spectral radius of $S$ is less
than one, then so is that of $A$, and $A$ is similar to a contraction.
\end{proposition}

\begin{proof}
Choose an integer $k$ with $\|S^k\|<1$.  By norm convergence,
$\|P_N^k\|<1$ for all sufficiently large $N$, so that $r(P_N)<1$ for
such $N$.  Fix one such $N\ge1$.  Since $P_N=R_N(A)$ for the rational
function $R_N=\prod_{j\le N}b_{a_j}$, whose poles lie outside
$\sigma(A)$, the spectral mapping theorem gives
\[
 \Bigl|\prod_{j=1}^N b_{a_j}(\lambda)\Bigr|<1
 \qquad(\lambda\in\sigma(A)).
\]
Each Blaschke factor has modulus one on $\T$ and modulus greater than
one outside $\overline\D$, wherever it is finite, because
$|z-a|^2-|1-\overline az|^2=(|z|^2-1)(1-|a|^2)$.  It follows that
$\sigma(A)\subset\D$, and by compactness $r(A)<1$.  Finally,
$Q=\sum_{n\ge0}A^{*n}A^n$ converges in norm, satisfies $Q\ge I$ and
$Q-A^*QA=I$, so $A$ is similar to a contraction.
\end{proof}

\subsection{Simultaneous transfer for commuting tuples}
\label{subsec:simultaneous}

\begin{theorem}[Simultaneous transfer]\label{thm:simultaneous}
Let $A_1,\dots,A_d$ be pairwise commuting bounded operators, let
$\theta_1,\dots,\theta_d$ be non-constant inner functions with
$\sigma(A_i)\cap\operatorname{Sp}(\theta_i)=\varnothing$, and put
$S_i=\theta_i(A_i)$.  If $S_1,\dots,S_d$ admit a common Lyapunov
operator $H$, that is, a coercive $H$ with
\begin{equation}\label{eq:common-metric}
 S_i^*HS_i\le H\qquad(1\le i\le d),
\end{equation}
then $A_1,\dots,A_d$ admit a common Lyapunov operator $Q$, that is, a
coercive $Q$ with $A_i^*QA_i\le Q$ for all $i$.
Consequently $(S_1,\dots,S_d)$ is simultaneously similar to a tuple of
contractions if and only if $(A_1,\dots,A_d)$ is.
\end{theorem}

\begin{proof}
Write $\theta_i=\eta_i B_iS_{\mu_i}$ and define a bounded positive
linear map on $\Bcal(\Hh)$ by
\begin{equation}\label{eq:Phi}
 \Phi_i(K)=Q_{B_i}(K)+B_i(A_i)^*Q_{\mu_i}(K)B_i(A_i),
\end{equation}
where $Q_{B_i}(K)$ is the series \eqref{eq:Q-series} and
$Q_{\mu_i}(K)$ the integral \eqref{eq:Qmu}, both formed with $A_i$ in
place of $A$; for an empty Blaschke product the first term is zero and
$B_i(A_i)=I$.  These series and integrals converge in operator norm for
every bounded $K$.  The computations in the proofs of
Lemma~\ref{lem:singular} and Theorem~\ref{thm:inner} give, for every
self-adjoint $K$,
\begin{equation}\label{eq:Phi-identity}
 \Phi_i(K)-A_i^*\Phi_i(K)A_i=K-S_i^*KS_i,
\end{equation}
and $\Phi_i(K)$ is coercive whenever $K$ is, because $\theta_i$ is
non-constant.  Every summand and every integrand in $\Phi_i(K)$ has the
form $V^*KV$ with $V$ a function of $A_i$, so an operator $R$ that
commutes with $A_i$ commutes with all these $V$, and therefore
\begin{equation}\label{eq:Phi-preserve}
 \Phi_i(K)-R^*\Phi_i(K)R=\Phi_i(K-R^*KR).
\end{equation}

Put $H_0=H$ and $H_i=\Phi_i(H_{i-1})$ for $1\le i\le d$.  We claim that
$H_i$ is coercive, that $S_k^*H_iS_k\le H_i$ for every $k$, and that
$A_k^*H_iA_k\le H_i$ for $k\le i$.  For $i=0$ this is the hypothesis
\eqref{eq:common-metric}.  At step $i$, the inequality for $A_i$
follows from \eqref{eq:Phi-identity} and the inductive hypothesis for
$S_i$, while the inequalities for $S_k$ (all $k$) and for $A_k$ ($k<i$)
are preserved by \eqref{eq:Phi-preserve} and the positivity of
$\Phi_i$, since these operators commute with $A_i$.  Hence $Q=H_d$ has
the required properties.  The converse follows from the converse part
of Theorem~\ref{thm:inner}, applied with one common similarity.
\end{proof}

%%=====================================================================
\section{Boundary obstructions}
\label{sec:negative}
%%=====================================================================

Theorem~\ref{thm:blaschke} is a statement about a norm-convergent
product under the separation hypothesis.  It must not be read as a
statement about an arbitrary $H^\infty$ functional calculus, and the
three examples of this section show why.  In each of them the
separation hypothesis fails at the point where the poles accumulate.

\subsection{A bounded calculus does not suffice}

\begin{proposition}\label{prop:pisier-gamal}
There exist a polynomially bounded operator $A$ on a Hilbert space that
is not similar to a contraction, and an infinite Blaschke product $B$
with zeros $a_j\in(0,1)$, such that $B(A)=0$.  The zeros, listed with
multiplicity, can moreover be chosen so that
$\sum_j\sqrt{1-a_j^2}<\infty$.  In this example
the point $1$ lies both in $\sigma(A)$ and in the closure of the set of
poles of $B$.
\end{proposition}

\begin{proof}
We use the modification of Pisier's counterexample due to Gamal'
\cite[\S5]{Gamal16}; see also the introduction of
\cite{GamalSufficient}.  By the finite-dimensional form of Pisier's
counterexample there are nilpotent operators $R_j$ on finite-dimensional
Hilbert spaces $H_j$ whose polynomial bounds are bounded by a common
constant $K$ while their complete polynomial bounds are unbounded along
a subsequence.  We write $R_j^{d_j}=0$.

For $0<w<1$ put $\beta_w(z)=(w-z)/(1-wz)$.  Disc automorphisms preserve
both the polynomial and the complete polynomial bounds, by rational
approximation on the closed disc in one direction and by
$\beta_w\circ\beta_w=z$ in the other.  Choose $w_j\uparrow1$ so rapidly
that
\[
 d_j\sqrt{1-w_j^2}\le2^{-j},\qquad
 \|\beta_{w_j}(R_j)-I_{H_j}\|\le2^{-j},
\]
the second condition being attainable because $I-R_j$ is invertible
and $\beta_w(R_j)\to I$ as $w\uparrow1$.  Set
\[
 A_j=\beta_{w_j}(R_j),\quad A=\bigoplus_jA_j,\quad
 B=\prod_j\beta_{w_j}^{d_j}.
\]
Then $B$ is a Blaschke product with the stated square-root summability,
and $A$ is polynomially bounded by $K$ while its complete polynomial
bound is infinite, so $A$ is not similar to a contraction by Paulsen's
theorem \cite{Pau}.  An $H^\infty$ calculus is defined by
$f(A)=\bigoplus_jf(A_j)$, where $f(A_j)$ is given by the Riesz--Dunford
calculus since $\sigma(A_j)=\{w_j\}\subset\D$.  The finite-dimensional
calculi satisfy $\|f(A_j)\|\le K\|f\|_\infty$: indeed $f_r(z)=f(rz)$ is
a uniform limit of polynomials on $\overline\D$ with
$\|f_r\|_\infty\le\|f\|_\infty$, and $f_r(A_j)\to f(A_j)$ in norm as
$r\uparrow1$.  The calculus is weak-$*$ continuous on bounded nets:
bounded weak-$*$ convergence in $H^\infty$ implies locally uniform
convergence in $\D$, hence norm convergence of $f_\alpha(A_j)$ for each
$j$, and a truncation of the direct sum then gives convergence in the
weak operator topology.  Since $\beta_{w_j}(A_j)=R_j$, the power
$\beta_{w_j}^{d_j}$ annihilates $A_j$; moreover
$B=\beta_{w_j}^{d_j}G_j$ with $G_j$ holomorphic near $w_j$, because the
remaining factors converge locally uniformly on $\D$, so
$B(A_j)=R_j^{d_j}G_j(A_j)=0$ for every $j$, that is, $B(A)=0$.
Finally, we have $\sigma(A)=\{w_j:j\ge1\}\cup\{1\}$, and the poles $w_j^{-1}$
also tend to $1$.
\end{proof}

Thus a bounded, weak-$*$ continuous $H^\infty$ calculus is compatible
with $B(A)=0$ and with the failure of similarity.  The refined
construction of \cite{Gamal16} even yields cyclic examples.  In the
other direction, for a Carleson--Newman Blaschke product $B$ an
absolutely continuous polynomially bounded operator $A$ with $B(A)=0$
\emph{is} similar to a contraction; see \cite[Theorem~3.2]{GamalSufficient}.

\subsection{A norm radial limit does not suffice}

\begin{example}\label{ex:radial}
Put $a_n=(n^2-1)/(n^2+1)$,
\[
 A=\bigoplus_{n\ge1}\begin{pmatrix}a_n&1\\0&a_n\end{pmatrix},
 \qquad
 F(z)=\prod_{n\ge1}\frac{a_n-z}{1-a_nz},\qquad B=F^2 .
\]
Then $A$ is bounded, $\sigma(A)=\{a_n:n\ge1\}\cup\{1\}$, and
$\|A^k\|\ge\sup_nka_n^{k-1}=k$, so $A$ is not power bounded.
Nevertheless $B(rA)\to0$ in operator norm as $r\uparrow1$, where each
$B(rA)$ is defined by the holomorphic functional calculus because
$\sigma(rA)$ is a compact subset of $\D$.
\end{example}

\begin{proof}
The product is a Blaschke product because
$\sum_n(1-a_n)=\sum_n2/(n^2+1)<\infty$.  Writing $w=(1+z)/(1-z)$, the
Euler product formulae give
\[
 F(z)=\prod_{n\ge1}\frac{1-w/n^2}{1+w/n^2}
 =\frac{\sin(\pi\sqrt w)}{\sinh(\pi\sqrt w)},
\]
the quotient being independent of the choice of the square root.  For
the $k$-th power,
$A^k=\bigoplus_n\begin{psmallmatrix}a_n^k&ka_n^{k-1}\\0&a_n^k
\end{psmallmatrix}$, so $\|A^k\|\ge\sup_nka_n^{k-1}=k$.

For real $x\uparrow1$ put $t=\sqrt{(1+x)/(1-x)}$.  Then
$F(x)=O(e^{-\pi t})$ and $F'(x)=O(t^3e^{-\pi t})$, so $B$ and $B'$
extend continuously to $[0,1]$ with $B(1)=B'(1)=0$; moreover
$B(a_n)=B'(a_n)=0$ for every $n$, the zeros being double.  Since
\[
 B(rA)=\bigoplus_{n\ge1}
 \begin{pmatrix}B(ra_n)&rB'(ra_n)\\0&B(ra_n)\end{pmatrix},
\]
the uniform continuity of $B$ and $B'$ on $[0,1]$ gives
$\|B(rA)\|\to0$.
\end{proof}

Thus the existence of a norm radial limit as $r\uparrow1$ does not
imply the power boundedness of $A$, even when that limit is zero.  The
example also shows that radial evaluation is quite different from the
norm product: here $1\in\sigma(A)$ is an accumulation point of the
poles $a_n^{-1}$, so the separation hypothesis fails and the norm
product is not available.

\begin{remark}\label{rem:radial-partial}
The partial products behave quite differently from $B(rA)$.  For
$B_N=\prod_{n\le N}((a_n-z)/(1-a_nz))^2$ and $n>N$ one has
$B_N(A)=B_N(a_n)I+B_N'(a_n)N_n$ on the $n$-th block, where $N_n$ is the
nilpotent part of that block, so $\|B_N(A)\|\ge|B_N'(a_n)|$; letting
$n\to\infty$ and using the continuity of $B_N'$ at $1$ gives
$\|B_N(A)\|\ge|B_N'(1)|=2\sum_{n\le N}n^2\to\infty$.  Norm radial
evaluation and the convergence of the partial products are therefore
unrelated hypotheses.
\end{remark}

\subsection{Strong products fail, even for Ritt operators}
\label{subsec:ritt-strong}

The third example shows that, even for Ritt operators, the strong
convergence of the partial products to zero is compatible with $A$ not
even being polynomially bounded.  By a theorem of Le Merdy
\cite[Theorem~5.1]{LeMerdy}, a polynomially bounded Ritt operator on a
Hilbert space is similar to a contraction; within the class of Ritt
operators, therefore, failure of similarity to a contraction is the same
as failure of polynomial boundedness, and the conclusion of the
following proposition is the strongest possible.  Ritt operators that
are not polynomially bounded were first constructed by Le Merdy
\cite[Proposition~5.2]{LeMerdy}, on the basis of an example of Baillon
and Cl\'ement; our operator $A$ is a variant of his construction, and
the multiplier formula used below is taken from the same proof.

\begin{proposition}\label{prop:ritt-strong}
There exist a Ritt operator $A$ on a Hilbert space that is not
polynomially bounded, and hence not similar to a contraction, and an
infinite Blaschke product $B$ whose zeros satisfy
$\sum_n\sqrt{1-|a_n|^2}<\infty$, such that the partial products
$B_N(A)$ converge strongly to zero.
\end{proposition}

\begin{proof}
Let $(x_n)$ be a conditional Schauder basis of a Hilbert space, with
partial-sum projections $P_n$ and $K=\sup_n\|P_n\|$; in this proof
$P_n$ denotes these projections rather than partial Blaschke products.
We shall use the following bounded-variation multiplier formula, which
is \cite[(5.3)]{LeMerdy}: if $c_n\to c_\infty$ and
$V(c)=\sum_n|c_{n+1}-c_n|<\infty$, then
\[
 M_c=c_\infty I-\sum_{n\ge1}(c_{n+1}-c_n)P_n,\qquad M_cx_n=c_nx_n,
\]
the series converging absolutely in norm, with
$\|M_c\|\le|c_\infty|+KV(c)$; indeed
$\sum_{n\ge m}(c_{n+1}-c_n)=c_\infty-c_m$.

Put $\lambda_n=1-2^{-n}$ and $Ax_n=\lambda_nx_n$.  The multiplier
formula shows that $A$ is bounded, and that
$\sigma(A)=\{\lambda_n\}\cup\{1\}$: for $\lambda$ outside this closed
set, the formula applied to $c_n=(\lambda-\lambda_n)^{-1}$, whose
variation is at most
$\sum_n(\lambda_{n+1}-\lambda_n)/\dist(\lambda,\{\lambda_n\}\cup\{1\})^2$,
produces the inverse of $\lambda I-A$.  Since $\lambda_n^k$ increases
to $1$, the formula also gives $\|A^k\|\le1+K$ for all $k$.  For
$g_k(t)=kt^{k-1}(1-t)$, which is non-negative, unimodal and bounded by
$1$ on $[0,1]$, the sequence $g_k(\lambda_n)$ is unimodal in $n$, tends
to $0$, and has variation at most $2$; hence
$k\|A^{k-1}(I-A)\|=\|M_{g_k}\|\le2K$ for all $k$.  By the
power-and-difference characterisation of Ritt operators
\cite{NZ,Lyu}, $A$ is a Ritt operator.

The sequence $(\lambda_n)$ is interpolating: for $|n-m|=j$,
\[
 \left|\frac{\lambda_n-\lambda_m}{1-\lambda_n\lambda_m}\right|
 \ge\frac{1-2^{-j}}{1+2^{-j}},
 \quad\text{so}\quad
 \inf_n\prod_{m\ne n}\left|\frac{\lambda_n-\lambda_m}
 {1-\lambda_n\lambda_m}\right|\ge\prod_{j\ge1}
 \left(\frac{1-2^{-j}}{1+2^{-j}}\right)^2>0 .
\]
By Carleson's interpolation theorem there is a constant $C$ such that
for every choice of signs $(\varepsilon_n)$ there is $f\in H^\infty$
with $\|f\|_\infty\le C$ and $f(\lambda_n)=\varepsilon_n$ for all $n$.

Suppose that $A$ were polynomially bounded, with bound $M$.  For
$0<r<1$ the function $f_r(z)=f(rz)$ is a uniform limit of polynomials
on the closed disc, so $\|f_r(A)\|\le MC$, and $f_r(A)$ converges
strongly, as $r\uparrow1$, to the sign multiplier
$x_n\mapsto\varepsilon_nx_n$, since it does so on finite linear
combinations of the basis vectors.  All sign multipliers would then be
bounded by $MC$, and the basis would be unconditional, a contradiction.
Hence $A$ is not polynomially bounded.

Finally, let $B=\prod_{n\ge1}b_{\lambda_n}^2$ and
$B_N=\prod_{j\le N}b_{\lambda_j}^2$.  For each $j\le N$ the function
$b_{\lambda_j}^2$ is non-negative and increasing on $[\lambda_N,1)$, so
for $n>N$ the sequence $c_n=B_N(\lambda_n)$ increases from $0$ to $1$,
while $c_n=0$ for $n\le N$; thus $c_\infty=1$, $V(c)\le1$, and
$\|B_N(A)\|\le1+K$ for all $N$.  Each $B_N(A)$ annihilates
$x_1,\dots,x_N$, so the uniform bound and the density of finite linear
combinations give $B_N(A)\to0$ strongly.  If the double zeros are
listed individually, the intermediate partial products also converge
strongly to zero: the same formula gives
$\|b_{\lambda_j}(A)\|\le1+2K$, because
$t\mapsto(\lambda_j-t)/(1-\lambda_jt)$ decreases along $(\lambda_n)$
from a value in $(-1,1)$ to the limit $-1$.
\end{proof}

Consequently the Ritt condition does not restore the implication
\eqref{eq:Q-intro} when an infinite Blaschke product is evaluated at
$A$ as a strong limit of partial products.  This is consistent with
Theorem~\ref{thm:blaschke}: here the poles $\lambda_n^{-1}$ accumulate
at $1\in\sigma(A)$.

%%=====================================================================
\section{Finite Blaschke products: polynomial Lyapunov operators}
\label{sec:finite}
%%=====================================================================

For a finite Blaschke product the Lyapunov operator can be written
without any resolvents, and this yields a quantitative bound that is
used in the next section.  Write
\begin{equation}\label{eq:ND}
 B(z)=\eta\frac{N(z)}{D(z)},\quad
 N(z)=\prod_{j=1}^m(z-a_j),\quad
 D(z)=\prod_{j=1}^m(1-\overline{a_j}z),\quad |\eta|=1,
\end{equation}
with $m\ge1$ and $a_j\in\D$.

\begin{proposition}[Polynomial Lyapunov operator]\label{prop:poly-metric}
Suppose that $D(A)$ is invertible and that $B(A)$ is similar to a
contraction with $\simc(B(A))\le K$.  Then there is a constant
$C_B<\infty$, depending only on $B$, such that
\begin{equation}\label{eq:finite-power}
 \simc(A)\le C_BK\max_{0\le k<m}\|A^k\|
 \le C_BK\max\{1,\|A\|\}^{m-1}.
\end{equation}
\end{proposition}

\begin{remark}
No hypothesis on $\sigma(A)$ beyond the invertibility of $D(A)$ is
needed in Proposition~\ref{prop:poly-metric}; in particular, the bound
\eqref{eq:finite-power} is uniform over all operators $A$ with a common
bound on $\|A\|$ and on $\simc(B(A))$.
\end{remark}

\begin{proof}[Proof of Proposition \ref{prop:poly-metric}]
By Proposition~\ref{prop:metric-p2}, there is a Lyapunov operator $H$
for $B(A)$ with $\sqrt{\|H\|\|H^{-1}\|}\le K$.  For
$j=1,\dots,m$ put
\begin{equation}\label{eq:Fj}
 F_j(z)=\sqrt{1-|a_j|^2}\prod_{\ell<j}(z-a_\ell)
 \prod_{\ell>j}(1-\overline{a_\ell}z)\in\mathcal P_{m-1},
\end{equation}
where $\mathcal P_{m-1}$ denotes the space of polynomials of degree at
most $m-1$, and write $n_\ell(z)=z-a_\ell$ and
$d_\ell(z)=1-\overline{a_\ell}z$.  The one-factor computation
\[
 \overline{d_\ell(w)}d_\ell(z)-\overline{n_\ell(w)}n_\ell(z)
 =(1-a_\ell\overline w)(1-\overline{a_\ell}z)
 -(\overline w-\overline{a_\ell})(z-a_\ell)=(1-\overline wz)(1-|a_\ell|^2)
\]
and the telescoping decomposition
\begin{multline*}
 \prod_{\ell=1}^m\overline{d_\ell(w)}d_\ell(z)
 -\prod_{\ell=1}^m\overline{n_\ell(w)}n_\ell(z)\\
 =\sum_{j=1}^m\Bigl(\prod_{\ell<j}\overline{n_\ell(w)}n_\ell(z)\Bigr)
 \bigl(\overline{d_j(w)}d_j(z)-\overline{n_j(w)}n_j(z)\bigr)
 \Bigl(\prod_{\ell>j}\overline{d_\ell(w)}d_\ell(z)\Bigr)
\end{multline*}
together give the Takenaka--Malmquist identity
\begin{equation}\label{eq:poly-kernel}
 \overline{D(w)}D(z)-\overline{N(w)}N(z)
 =(1-\overline wz)\sum_{j=1}^m\overline{F_j(w)}F_j(z),
\end{equation}
an identity between polynomials in $z$ and $\overline w$.

Replacing $z$ by $A$ and $\overline w$ by $A^*$, and inserting $H$
between the adjoint and the non-adjoint factors, turns
\eqref{eq:poly-kernel} into
\begin{equation}\label{eq:poly-stein}
 \widehat Q-A^*\widehat QA
 =D(A)^*HD(A)-N(A)^*HN(A)
 =D(A)^*\bigl[H-B(A)^*HB(A)\bigr]D(A)\ge0,
\end{equation}
where $\widehat Q=\sum_{j=1}^mF_j(A)^*HF_j(A)$ and we have used
$N(A)=\overline\eta B(A)D(A)$.

The polynomials $F_1,\dots,F_m$ form a basis of $\mathcal P_{m-1}$.
Indeed, if $\sum_jc_jF_j=0$, then evaluation at $a_1$ gives $c_1=0$,
since $F_1(a_1)\ne0$ while every other summand contains the factor
$z-a_1$; dividing by $z-a_1$ and evaluating at $a_2$ gives $c_2=0$, and
so on, repeated zeros causing no difficulty.  Hence there are scalars
$\alpha_j$ with $\sum_j\alpha_jF_j=1$, and with
$c_\alpha=(\sum_j|\alpha_j|^2)^{1/2}$ we have
\[
 \|x\|^2\le c_\alpha^2\sum_j\|F_j(A)x\|^2
 \quad\text{and therefore}\quad
 \widehat Q\ge\frac{1}{c_\alpha^2\|H^{-1}\|}\,I .
\]
Writing $F_j(z)=\sum_{k<m}c_{jk}z^k$,
$c_F=\bigl(\sum_j(\sum_k|c_{jk}|)^2\bigr)^{1/2}$ and
$M_A=\max_{k<m}\|A^k\|$, we also have
$\widehat Q\le\|H\|c_F^2M_A^2I$.  Thus $\widehat Q^{1/2}A\widehat Q^{-1/2}$
is a contraction, with condition number at most $c_\alpha c_FKM_A$,
which is \eqref{eq:finite-power} with $C_B=c_\alpha c_F$.
\end{proof}

\begin{remark}[Sharpness of the exponent]\label{rem:exponent}
For $B(z)=z^m$ with $m\ge2$, let $J_m$ be the nilpotent shift on
$\C^m$ and $A_t=tJ_m$ with $t\ge1$.  Then $B(A_t)=0$, so
$\simc(B(A_t))=1$, while any similarity of $A_t$ to a contraction has
condition number at least $\|A_t^{m-1}\|=t^{m-1}$.  Hence the exponent
$m-1$ in a bound in terms of $\|A\|$ alone cannot be lowered.  For
$B(z)=z^m$ the construction above reproduces the elementary Lyapunov
operator $\sum_{k<m}A^{*k}HA^k$.
\end{remark}

%%=====================================================================
\section{Uniform strong test collections}
\label{sec:uniform}
%%=====================================================================

This section is inspired by the work of Dritschel, Est\'evez and
Yakubovich \cite{DEY}.  Fix non-constant finite Blaschke products
$B_1,\dots,B_n$, write
$\Psi=(B_1,\dots,B_n)$, and let $\mathcal A_\Psi$ be the algebra of
finite sums of functions of the form
\begin{equation}\label{eq:A-Psi}
 \prod_{j=1}^nf_j(B_j(z)),\qquad f_j\in A(\overline\D),
\end{equation}
where $A(\overline\D)$ is the disc algebra.  We say that $\Psi$ is a
\emph{uniform strong test collection} on $\overline\D$ if for every
$K\ge1$ there is a constant $L(\Psi,K)$ such that
\[
 \sigma(A)\subset\overline\D,\quad
 \simc(B_j(A))\le K\ \text{ for all } j
 \quad\Longrightarrow\quad
 \simc(A)\le L(\Psi,K),
\]
where the similarities for different $j$ are not required to coincide;
and that $\Psi$ is a \emph{uniform test collection} if the same holds
with $\|B_j(A)\|\le1$ in place of $\simc(B_j(A))\le K$.

\begin{theorem}[Uniform strong testing from a dense algebra]
\label{thm:dense}
Suppose that $\mathcal A_\Psi$ is dense in $A(\overline\D)$ in the
uniform norm.  Then $\Psi$ is a uniform strong test collection on
$\overline\D$.
\end{theorem}

\begin{proof}
Let $m=\deg B_1$, let $C=C_{B_1}$ be the constant of
Proposition~\ref{prop:poly-metric}, and put
$M_A=\max_{k<m}\|A^k\|$.  That proposition shows that $A$ is similar to
a contraction with $\simc(A)\le CKM_A$, so that
$\|g(A)\|\le CKM_A\|g\|_\infty$ for every $g\in A(\overline\D)$; this
bound depends on $A$ and will be used only to control approximation
errors.

Choose $\eps=(2CK)^{-1}$ and, by density, $g_k\in\mathcal A_\Psi$ with
$\|z^k-g_k\|_\infty<\eps$ for $0\le k<m$, where $g_0=1$.  Writing
$g_k=\sum_{\ell=1}^{r_k}\prod_{j=1}^nf_{k\ell j}(B_j(z))$, we have
\[
 \|g_k(A)\|\le\sum_{\ell=1}^{r_k}K^n\prod_{j=1}^n
 \|f_{k\ell j}\|_\infty\le G,
\]
where $G$ is the largest of these finitely many bounds over
$0\le k<m$; here we have used that each $f_{k\ell j}(B_j(A))$ has norm
at most $K\|f_{k\ell j}\|_\infty$, by the similarity of $B_j(A)$ to a
contraction and von Neumann's inequality.  The bound $G$ depends on
$\Psi$, $K$ and the chosen $g_k$, but not on $A$.  On the other hand,
$\|A^k-g_k(A)\|\le CKM_A\eps=M_A/2$, so $M_A\le G+M_A/2$, that is,
$M_A\le2G$.  Applying Proposition~\ref{prop:poly-metric} once more
gives $\simc(A)\le2CKG$, which is the required uniform bound.
\end{proof}

\begin{corollary}[Injective unramified tuples]\label{cor:injective}
Let $\Psi=(B_1,\dots,B_n)$ be a tuple of non-constant finite Blaschke
products that is injective on $\overline\D$ and satisfies
\begin{equation}\label{eq:joint-derivative}
 (B_1'(z),\dots,B_n'(z))\ne0\qquad(z\in\overline\D).
\end{equation}
Then $\Psi$ is a uniform strong test collection on $\overline\D$.
\end{corollary}

\begin{proof}
We verify the density hypothesis of Theorem~\ref{thm:dense} by means of
the algebra-generation results of \cite{DEY}.  Finite Blaschke products
extend holomorphically across $\T$ with non-vanishing derivative there,
so the injectivity of $\Psi$ on $\T$ makes $\Psi$ weakly admissible in
the sense of \cite{DEY}.  By \cite[Lemma~17]{DEY} there is an admissible
family $\Phi$ whose coordinates have the form
$\varphi_\ell(z)=\prod_{j=1}^nh_{\ell j}(B_j(z))$ with
$h_{\ell j}\in A(\overline\D)$ and $\|h_{\ell j}\|_\infty\le1$.  Adjoin
the original coordinates to $\Phi$, assigning empty boundary arcs to the
added functions, as in the proof of \cite[Theorem~18]{DEY}.  The
enlarged family $\Phi_0$ is admissible, injective, and has
non-vanishing joint derivative, so \cite[Theorem~II]{DEY} gives
$\mathcal A_{\Phi_0}=A(\overline\D)$.

Let $\mathcal C$ be the uniform closure of $\mathcal A_\Psi$.  It
contains every coordinate of $\Phi_0$; and if $\varphi$ is such a
coordinate and $f\in A(\overline\D)$, then uniform polynomial
approximation of $f$ shows that $f\circ\varphi\in\mathcal C$.  Hence
$\mathcal A_{\Phi_0}\subset\mathcal C$, so $\mathcal C=A(\overline\D)$,
and Theorem~\ref{thm:dense} applies.
\end{proof}

\begin{theorem}[Exact characterisation]\label{thm:characterisation}
For a tuple $\Psi=(B_1,\dots,B_n)$ of non-constant finite Blaschke
products the following are equivalent:
\begin{enumerate}[label=\textup{(\roman*)}]
\item $\Psi$ is injective on $\overline\D$ and
      \eqref{eq:joint-derivative} holds;
\item $\Psi$ is a uniform strong test collection on $\overline\D$;
\item $\Psi$ is a uniform test collection on $\overline\D$.
\end{enumerate}
\end{theorem}

\begin{proof}
The implication (i)$\Rightarrow$(ii) is Corollary~\ref{cor:injective},
and (ii)$\Rightarrow$(iii) is immediate, since $\|B_j(A)\|\le1$ is the
case $K=1$ of $\simc(B_j(A))\le K$.

For (iii)$\Rightarrow$(i), suppose first that $z\ne w$ in $\overline\D$
satisfy $B_j(z)=B_j(w)$ for all $j$, and for $t>0$ put
\[
 A_t=\begin{pmatrix}z&t(w-z)\\0&w\end{pmatrix}.
\]
Then $A_t$ is diagonalisable with eigenvalues $z\ne w$, so
$B_j(A_t)=B_j(z)I$ for every $j$, a contraction, while
$\|A_t\|\ge t|w-z|\to\infty$.  A uniform test collection would give
$\simc(A_t)\le L(\Psi,1)$, whereas $\simc(A_t)\ge\|A_t\|$; this is
impossible.  If instead the joint derivative vanishes at some
$z_0\in\overline\D$, take $A_t=z_0I+tN$ with
$N=\begin{psmallmatrix}0&1\\0&0\end{psmallmatrix}$; then
$B_j(A_t)=B_j(z_0)I+tB_j'(z_0)N=B_j(z_0)I$, again a contraction, while
$\|A_t\|\ge t\to\infty$.
\end{proof}

\begin{corollary}[Common peripheral collapse]\label{cor:collapse}
Let $E\subset\T$ be a finite set with at least two points, and suppose
that $B_j(\zeta)=1$ for every $\zeta\in E$ and every $j$.  Then $\Psi$
is not a uniform test collection on $\overline\D$, even when attention
is restricted to $\mathrm{Ritt}_E$ operators, with no bound imposed on
their $\mathrm{Ritt}_E$ constants.
\end{corollary}

\begin{proof}
Take $z=\zeta_1\ne\zeta_2=w$ in $E$ in the first matrix of the
preceding proof.  Then $B_j(A_t)=I$ for all $j$, and $A_t$ is similar
to $\operatorname{diag}(\zeta_1,\zeta_2)$:
\[
 (\lambda I-A_t)^{-1}=P_t\operatorname{diag}
 \bigl((\lambda-\zeta_1)^{-1},(\lambda-\zeta_2)^{-1}\bigr)P_t^{-1},
 \qquad P_t=\begin{pmatrix}1&t\\0&1\end{pmatrix}.
\]
Hence $\dist(\lambda,E)\|(\lambda I-A_t)^{-1}\|\le\|P_t\|\|P_t^{-1}\|$
for $|\lambda|>1$, so $A_t$ is a $\mathrm{Ritt}_E$ operator, whereas
$\simc(A_t)\ge\|A_t\|\ge t|\zeta_1-\zeta_2|$ is unbounded.
\end{proof}

\begin{example}[Three tests, none injective]\label{ex:triple}
Fix $0<a<1$ and put $B_1(z)=z^2$, $B_2(z)=z^3$ and
$B_3(z)=((z-a)/(1-az))^2$.  The pair $(z^2,z^3)$ is injective on
$\overline\D$, and its joint derivative vanishes only at $0$; since
$B_3'(0)=-2a(1-a^2)\ne0$, the joint derivative of the triple vanishes
nowhere.  Thus Corollary~\ref{cor:injective} applies, although no
single test is injective.

In this example the constant can be computed.  With $u=z^2$, $v=z^3$
and $w=B_3(z)$,
\begin{equation}\label{eq:triple-formula}
 z=v+\frac{a^2(u-1)^2+(1-a^2)^2u-w(1-a^2u)^2}{2a(1-a^2)} .
\end{equation}
Indeed, $w(1-a^2z^2)^2=\bigl[a(z^2-1)+(1-a^2)z\bigr]^2$, so the
numerator equals $-2a(1-a^2)z^3+2a(1-a^2)z$, and
\eqref{eq:triple-formula} is an identity.  Suppose that
$\simc(B_j(A))\le K$ for $j=1,2,3$.  Then $\|A^3\|\le K$ and
$\|B_3(A)\|\le K$; moreover, with $q(u)=a^2(u-1)^2+(1-a^2)^2u$, the
similarity of $A^2$ to a contraction and von Neumann's inequality give
$\|q(A^2)\|\le K\|q\|_\infty\le K(1+a^2)^2$ and
$\|(I-a^2A^2)^2\|\le K(1+a^2)^2$.  Evaluating
\eqref{eq:triple-formula} at $A$ therefore gives
\[
 \|A\|\le L(a,K):=K+\frac{K(1+K)(1+a^2)^2}{2a(1-a^2)} .
\]
If $H$ is a Lyapunov operator for $A^2$ with
$\sqrt{\|H\|\|H^{-1}\|}\le K$, which exists by
Proposition~\ref{prop:metric-p2}, then $Q=H+A^*HA$ satisfies
$Q-A^*QA=H-(A^2)^*HA^2\ge0$ and $H\le Q\le(1+\|A\|^2)\|H\|I$, so that
\begin{equation}\label{eq:triple-bound}
 \simc(A)\le K\sqrt{1+L(a,K)^2}.
\end{equation}
\end{example}

\begin{remark}[What is used and what is not]
In Example~\ref{ex:triple} the bound $\|A^3\|\le K$ comes from the
similarity of $B_2(A)=A^3$, and $\|B_3(A)\|\le K$ from that of
$B_3(A)$; neither follows from the similarity of $A^2$ alone.  The
final Lyapunov operator in \eqref{eq:triple-bound}, however, uses only the
similarity of $A^2$ together with the crude norm bound on $A$, which is
why the constant is so large.
\end{remark}

\begin{remark}[Relation to \cite{DEY}]
Corollary~19 of \cite{DEY} gives a uniform \emph{test} collection for
injective unramified families of inner functions on finitely connected
analytic domains, and the discussion on p.~1017 of \cite{DEY} leaves
the strong version open.  Theorem~\ref{thm:characterisation} settles
the strong version on the disc, by an argument that is special to the
disc: a single non-constant finite Blaschke coordinate already makes
$A$ similar to a contraction, with a bound involving finitely many
powers of $A$, and the absorption argument in the proof of
Theorem~\ref{thm:dense} upgrades this to a uniform bound. 
\end{remark}

%%=====================================================================
\section{Declaration of AI use and acknowledgements}
\label{sec:AI}
%%=====================================================================

All mathematical content and arguments were developed solely by the author. 
Large Language Models assistance was limited strictly to the writing and proofreading stages of this manuscript. 
This work was supported by the CDP C2EMPI and its institutional partners through
the project R-CDP-24-004-C2EMPI and by the EU COST network OC-2025-1-28418.

%%=====================================================================

\end{document}